\documentclass[11pt,twoside]{article}
\usepackage{mathbbold}
\usepackage{amsmath, amssymb, mathrsfs, graphicx}

\def\scr{\mathscr}

\def\az{\alpha}  \def\bz{\beta}
    \def\dz{\delta}
\def\ez{\eta}    \def\fz{\varphi}
\def\gz{\gamma}  
\def\lz{\lambda}

        \def\uz{\theta}
\def\vz{\varepsilon}

\def\qd{\quad}
\def\qqd{\qquad}

\newcommand{\mathsym}[1]{{}}

\def\scr{\mathscr}

\def\le{\leqslant}
\def\ge{\geqslant}

\def\geqs{\geqslant}

\font\cms=cmss9 scaled \magstep1

\newcommand{\nnb}{\nonumber}
\def\nnd{\noindent}

\def\thm{\nnd\bg{thm1}}

\def\lmm{\nnd\bg{lmm1}}
\def\prp{\nnd\bg{prp1}}

\def\rmk{\nnd\bg{rmk1}}

\def\dethm{\end{thm1}}

\def\delmm{\end{lmm1}}
\def\deprp{\end{prp1}}

\def\dermk{\end{rmk1}}

\def\prf{\medskip \noindent {\bf Proof}. }
\def\deprf{\quad $\square$ \medskip}
\def\bg{\begin}
\def\be{\bg{equation}}
\def\de{\end{equation}}

\def\dear{\end{eqnarray}}
\def\lb{\label}
\def\ct{\cite}

\newcommand{\rf}[2]{[\ref{#1}; #2]}

\def\den{\end{enumerate}}

\def\d{\text{\rm d}}
\def\pp{\partial}
\def\<{\langle} \def\>{\rangle}

\def\Ric{\text{\rm Ric}}

\begin{document}

\allowdisplaybreaks[4]
\thispagestyle{empty}
\renewcommand{\thefootnote}{\fnsymbol{footnote}}

\noindent {Front. Math. China 2011, 6(6): 1025-1043}

\vspace*{.5in}
\begin{center}
{\bf\Large General estimate of the first eigenvalue on manifolds}
\vskip.15in {Mu-Fa Chen}
\end{center}
\begin{center} (Beijing Normal University, Beijing 100875, China) \end{center}
\vskip.1in

\begin{center}{December 30, 2010}\end{center}

\markboth{\sc Mu-Fa Chen}{\sc General estimate of the first eigenvalue on manifolds}


\footnotetext{{\it Key words and phases}.\quad First non-trivial eigenvalue,
sharp estimate, Riemannian manifold.}
\footnotetext{Research supported in part by the
         National Natural Science Foundation of China (Nos. 10721091, 11131003), by the ``985'' project from the Ministry
of Education in China.}

\bigskip

\begin{abstract}
Ten sharp lower estimates of the first non-trivial eigenvalue of Laplacian on compact Riemannian manifolds are reviewed and compared. An improved variational formula,
a general common estimate, and a new sharp one are added. The best lower estimates are now updated. The new estimates provide a global picture of what one can expect by our approach.\end{abstract}


\section{Introduction}


Let $M$ be a compact, connected Riemannian manifold, without or
with convex boundary $\pp M$. When $\pp M\ne\emptyset$, we adopt
Neumann boundary condition. Next, let $\Ric_M\geqs K$ for some
$K\in \mathbb R$. Denote by $d$ and $ D$, respectively, the dimension
and diameter of $M$. We are interested in the estimate of  the
first non-trivial eigenvalue $\lz_1$ of Laplacian. On this topic, there is
a great deal of publications (see for instance \ct{cmf05a},
Schoen and Yau \ct{srys},
Wang \ct{wfy04}, and references therein).
Throughout this paper, we use the quantity
$$\az=\az(K, d, D)=\frac{D}{2}\sqrt{\frac{-K}{d-1}},$$
which involves all of the three geometric quantities: $d$, $D$ and $K$.
Clearly, $|\az|=\az (|K|, d, D)$.
By the Myers theorem, we have $|\az|\le \pi/2$ whenever $K>0$ for a complete
Riemannian manifold. Therefore, the quantity $\az$ is geometrically meaningful.
We adopt the convention: $\az=0$ if $d=1$. With the necessary notation in mind,
the main result Theorem \ref{t13} and its illustrating figures 7--9 should be
readable now. One may have a look before going further.

We are going to recall the sharp lower estimates, only the related part of them for saving space.
First, for nonnegative curvature, the following sharp lower bounds are perhaps well known.
\begin{itemize}
\item Lichnerowicz (1958):
\be\dfrac{d}{d-1}\, K=\dfrac{4 d}{D^2}\, |\az|^2, \qquad d>1,\; K\ge 0. \lb{01}\de
\item B\'erard, Besson and Gallot (1985):
 \be d\bigg\{\dfrac{\int_0^{\pi/2}{\mbox{cos}}^{d-1}t \mbox{d} t}{ \int_0^{D/2}{\mbox{cos}}^{d-1}t \mbox{d} t}\bigg\}^{2/d}
 =d\bigg\{\dfrac{\int_0^{\pi/2}{\mbox{cos}}^{d-1}t \mbox{d} t}{ \int_0^{|\az|}{\mbox{cos}}^{d-1}t \mbox{d} t}\bigg\}^{2/d}, \qqd K=d-1>0.  \lb{02}\de
 Here and in what follows, $\cos^k t=(\cos t)^k$.
\item Chen and Wang (1997):
\be \dfrac{dK}{(d-1) (1-\cos^d |\az|)}
=\dfrac{4d |\az|^2}{D^2 (1-\cos^d |\az|)}, \qd d>1, \; K\ge  0. \lb{03}\de
\item Zhong and Yang (1984):
\be {\pi^2}\big/{D^2}, \qquad K\ge 0.\lb{04}
 \de
 \end{itemize}
It is clear that (\ref{02}) improves (\ref{01}) by the Myers theorem. Even though it is not so obvious but it is true that (\ref{03}) improves (\ref{02}).
The first three results indicate a long period for the improvements
of (\ref{01}) step by step. All of them are
sharp for the unit sphere in two or higher dimensions but fail for the unit circle $(K=0)$. To which, the sharp estimate is given by the last result (\ref{04}).
Note that when $|\az|\downarrow 0$, the limit of (\ref{03}) equals $8D^{-2}\in (0, \pi^2 D^{-2})$.

Secondly, consider the non-positive curvature in which case,
the problem becomes harder. Here are the main known sharp lower bounds.
\begin{itemize}
\item Yang (1990), Jia (1991), Chen and Wang (1994):
\be \dfrac{\pi^2}{D^2}\, e^{-(d-1)\az}, \qqd K\le 0. \lb{05} \de
\item Chen and Wang (1997):
\be \frac{1}{D^2}\sqrt{\pi^4+8(d-1)\az^2}\,\cosh^{1-d} \az, \qquad d>1, \; K\le 0. \lb{06} \de
\item \rf{cmf94}{1994, Theorem 6.6}(corrected version):
\be \frac{1}{D^2}\big({(d-1)\az} \tanh \az\;
\mbox{sech}\, \uz\big)^2, \qquad d>1, \; K\le 0,\lb{07}
 \de
 where $\uz$ is obtained in the following way. Let
$$\uz_1=2^{-1}{(d-1)\az}\,\tanh \az,\qqd
\uz_n=\uz_1\tanh \uz_{n-1},\qqd n\ge 2,$$
then $\uz_n\downarrow \uz$.
\end{itemize}
The first two results have the same decay rate but (6) is better than (5) in general. They are sharp for the unit circle but not the last result which is designed for large $\az$. The last two estimates are not comparable. For fixed $d$, (\ref{06}) is better than (\ref{07}) for smaller
$\az$ but the inverse assertion happens for large $\az$ (cf. Fig. 8 below).




Thirdly, consider the optimal linear approximation of the lower estimates with respect
to the curvature $K$. In other words, one looks for a good combination of the optimal
estimates (\ref{01}) and (\ref{04}).
Many authors have contributed to the result (\ref{08}) below. It was proved by
Zhao (1999) under the restriction $-5 \pi^2/(3D^2)\le K\le 0$ and the remainder gap was covered by Xu and Pang (2001) in the case of $K\le 0$. The case of $K>0$ was proved by Xu, Yang and Xu (2002). Independently, the assertion was proved, with computer assisted, for all real $K$
in Chen, Scacciatelli and Yao (2001) (where some refined estimates are included).
A more direct analytic proof with some improvement was given by Shi and Zhang (2007).
Recently, the result (\ref{08}) below has been reproved by Ling (2006, 2007) using a different approach.
\bg{itemize}
\item  The following lower bound is studied/proved in the papers just mentioned:
\be {\pi^2}\big/{D^2}+ K/2, \qqd K\in {\mathbb R}.  \lb{08}\de
\item More precisely, the lower bound given by Shi and Zhang (2007) is as follows:
\begin{align}
&\sup_{s\in (0, 1)} s\bigg[4 (1-s){\pi^2}\big/{D^2}+ K\bigg]\nnb\\
&\qd ={
\begin{cases}
\displaystyle\bigg(\dfrac{\pi}{D}+ \dfrac{K D}{4 \pi}\bigg)^2,\qd & -4\pi^2\le KD^2\le 4\pi^2\\
K, & KD^2\in\big( 4\pi^2, (d-1) \pi^2\big]\\
0, & KD^2<- 4\pi^2
\end{cases}
\lb{09}}\end{align}
\item A refined lower bound given by Chen, Scacciatelli and Yao (2001) is the following:
\be \frac{\pi^2}{D^2}+ \frac{K}{2} +\big(10-{\pi^2}\big)
\frac{K^2 D^2}{16},\qqd |K| D^2\le 4.\lb{10}\de
\end{itemize}
To see that (\ref{09}) improves (\ref{08}), simply set $s=1/2$. On the region
$$\{(K, D): |K| D^2\le 4\},$$
it is obvious that (\ref{10}) is better than (\ref{09}).
Unlike the results (\ref{01})--(\ref{07}), the bounds given in (\ref{08})---(\ref{10}) are independent of the dimension $d$.

Before moving further, let us make some remarks on (\ref{10}). Since one is seeking for the dimension-free estimate and
$$(d-1)\az\,\tanh(\az r)\; \big\uparrow\; -{K D^2} r/4=: -2\bz r\qqd\text{as } d\uparrow \infty$$
\begin{center}{{\includegraphics{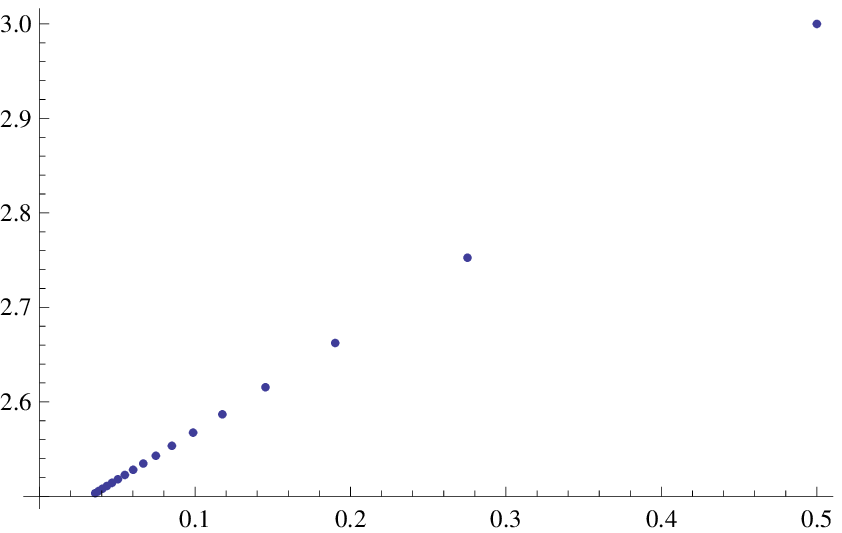}\hskip-3.7truecm
}
\vskip-3.0truecm{\hskip-0.2truecm\includegraphics{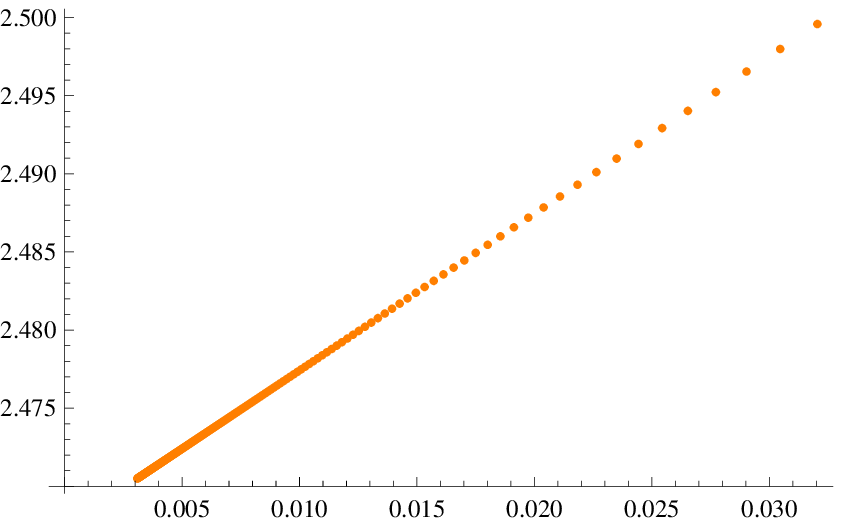}
}
\vskip-2.8truecm{\hskip-1.8truecm\includegraphics{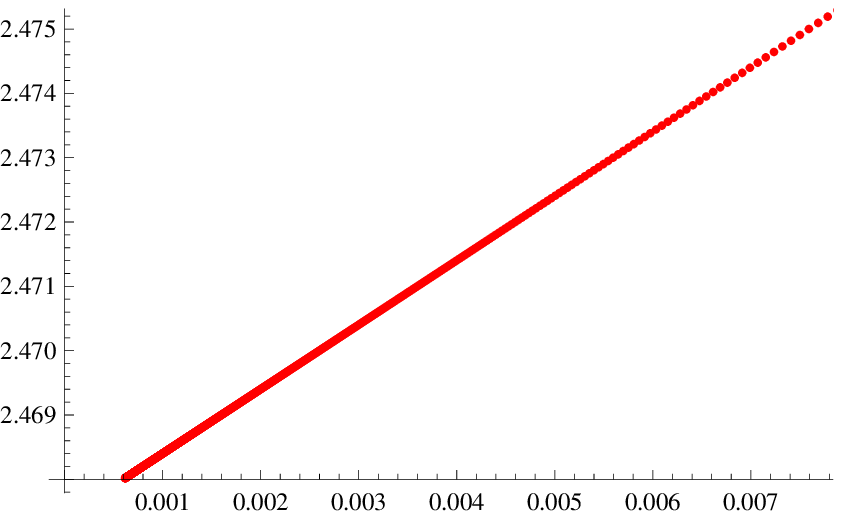}
}\newline
{\bf Figures 1--3}\qd\rm The first eigenvalue $\lz_0$ corresponding to $\{\bz_n\}_{n=2}^{18}$ (the top curve),
$\{\bz_n\}_{n=20}^{200}$ (the middle curve),
and $\{\bz_n\}_{n=2}^{1000}$ (the bottom curve), respectively.}\end{center}
(cf. proof of Theorem \ref{t17} below), the study on $\lz_1$ can be reduced to study the first eigenvalue (say $\lz_0$, for a moment) of the operator $\d^2/\d r^2-2\bz r \d /\d r$ on $(0, 1)$ (cf. \rf{csy}{Lemma 2.3} in which the constant $\az$ is replaced by $\bz$ here). This now  becomes a one-parameter problem and the required estimate (\ref{10}) becomes a product of ${4}{D^{-2}}$ and
\be{\pi^2}/{4} + \beta +(10-\pi^2) \beta^2,\qqd |\beta|\le  1 /2\lb{20}\de
which estimates $\lz_0$ from below. For a
sequence of $\bz$, say $\{\bz_n\}_{n=2}^{\infty}$ with $\bz_2=1/2$ and $\bz_n\downarrow 0$ as $n\to\infty$, the problem is numerically solvable and moreover, we have analytic solutions for the first four of $\{\bz_n\}$.
Actually, the corresponding eigenfunctions are all polynomials  (cf. \rf{csy}{Lemma 2.3}.
Now, (\ref{10}) is designed to be exact at $\bz=0$ and $\bz=1/2$
with the coefficient of the term $K$ to be $1/2$. The first eigenvalue $\lz_0$ just mentioned corresponding to $\{\bz_n\}$ is shown by three pictures for different region
of $n$ (Figures 1--3), and then the difference between the eigenvalue $\lz_0$ and its lower estimate (\ref{20}) is shown by two pictures (Figures 4 and 5).
In Fig.3, even though the whole sequence $\{\bz_n\}_{n=2}^{1000}$ is computed, but the output is restricted to a smaller interval near zero. To see the other part of the interval, one needs two more figures (1 and 2).
\medskip

\begin{center}{{{\includegraphics{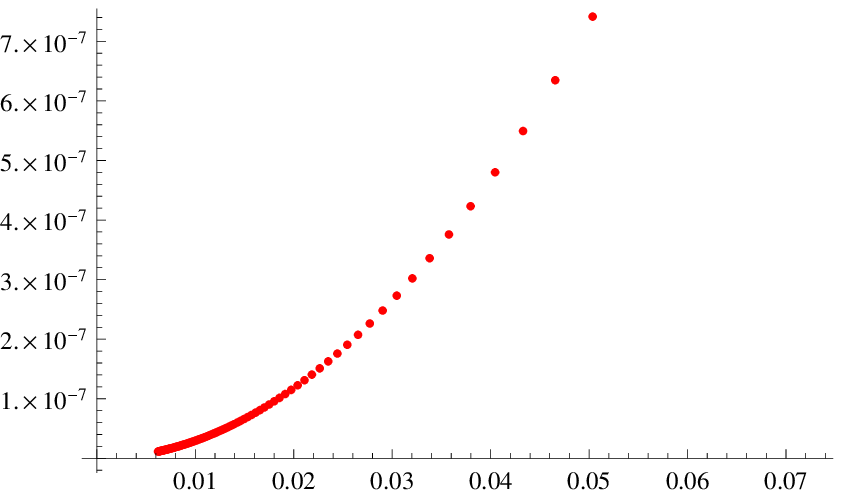}\hskip-3.8truecm}
}
\vskip-3.7truecm{\hskip-1.8truecm\includegraphics{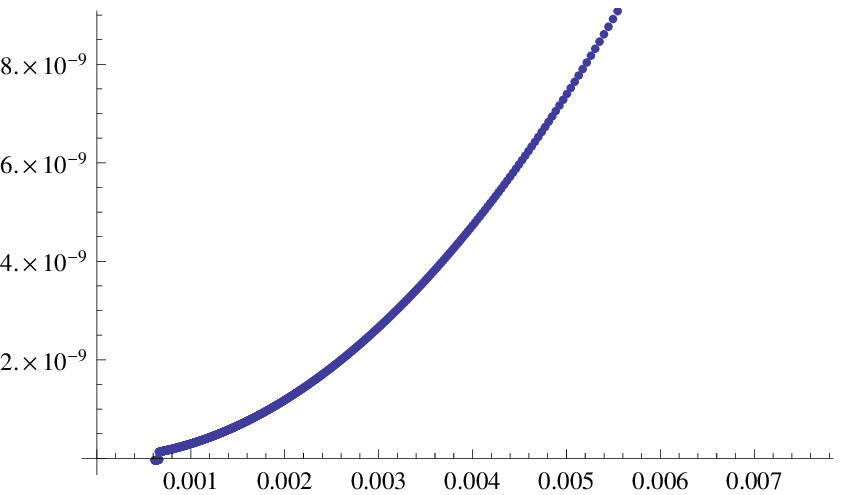}}\newline
{\bf Figure 4--5}\qd\rm The difference of the first eigenvalue $\lz_0$ and its lower estimate (\ref{20}) corresponding to $\{\bz_n\}_{n=2}^{30}$ (the curve on right) and
$\{\bz_n\}_{n=2}^{1000}$ (the curve on left), respectively.}\end{center}

The difficulty is that for each $n\ge 2$, one has to find a root of a polynomial having order $n-1$.
Here, instead of an analytic proof, we have used five figures to show carefully that the lower bound (\ref{10}) is rather sharp and the coefficient $1/2$ for the linear approximation should be exact. Besides, these figures also indicate that the analytic proofs would be too heavy for the paper and may not be essential, this is the reason why we often use figures, here and in what follows.
We mention that by a simple transform, the same conclusions hold for the sequence $\{-\bz_n\}_{n=1}^{\infty}$ (cf. \rf{csy}{Lemma 2.4}).
Therefore, (\ref{10}) holds for all real $K$ with $|K| D^2\le 4$.

Except the results (\ref{01})---(\ref{03}) and (\ref{07}), all of the above results are sharp only at one point ($K=0$), some of them can be very poor in some region.
For instance, (\ref{01})---(\ref{03}) are meaningless when $K=0$ and
(\ref{08})---(\ref{10}) can be negative or zero for sufficiently large $-K$.
The question now is the existence of a universal estimate. Fortunately,
the answer is affirmative as shown by Proposition \ref{t11} below.
Its first assertion implies the lower bounds (1)--(10), in terms of $\bar\lz$.

Recall that even though $\az$ is an imaginary number when $K>0$, the quantities
$\tanh(\az r)$ and $\cosh(\az r)$ used below are still meaningful: $\cosh(i\uz)= \cos \uz$, $\tanh(i\uz)= i\tan \uz$ for real $\uz$.

\prp\lb{t11} {\cms We have
$\lz_1\ge {4}{\bar\lz}/ {D^2}$ with the following estimates:
$$\frac{1}{4\dz}\le \frac{1}{\dz_1\wedge \dz_1^*}\le {\bar\lz}\le \frac{1}{\dz_1'\vee {\dz_1^*}'}\le \frac{1}{\dz},$$
where
$$\aligned
\dz&= \sup_{r\in (0, 1)} [\fz \psi](r),\\
\dz_1&=\sup_{r\in (0, 1)}\bigg\{\frac{1}{\sqrt{\fz (r)}}\int_0^r C \fz^{3/2}+
  \sqrt{\fz (r)}\int_r^1 C \fz^{1/2} \bigg\},\\
\dz_1'&=\sup_{r\in (0, 1)}\bigg\{\frac{1}{\fz (r)}\int_0^r C \fz ^{2} +
 [\fz  \psi] (r) \bigg\},\\
\dz_1^*&=\sup_{r\in (0, 1)}\bigg\{\frac{1}{\sqrt{\psi (r)}}\int_r^1 C^{-1} \psi^{3/2}+
  \sqrt{\psi (r)}\int_0^r C^{-1}\psi^{1/2} \bigg\},\\
{\dz_1^*}'&=\sup_{r\in (0, 1)}\bigg\{\frac{1}{\psi (r)}\int_r^1 C^{-1}\psi^{2}+
 [\fz \psi] (r)  \bigg\}
 \endaligned$$
(here the Lebesgue measure ``$\d u$'' is omitted) with
$$C(s)=\cosh^{d-1} (\az s),\qqd \fz (r)=\int_0^r C(u)^{-1}\d u,\qqd \psi (r)=\int_r^1 C(u)\d u.$$
All of the quantities used here depend on $d$ and $\az$.
}\deprp

To illustrate the power of Proposition \ref{t11}, consider the simplest case that $K=0$. Then $\dz=1/4$,
$\dz_1=\dz_1^*=5^{1/3}/4\approx 0.427$, $\dz_1'={\dz_1^*}'=3/8$, and so
 $$\frac{\dz_1}{\dz_1'}=\frac{\dz_1^*}{{\dz_1^*}'}=\frac{5^{1/3}}{4}\bigg /\frac{3}{8}\approx 1.14.$$
 The sharp estimate for $\lz_1$ is $\pi^2/D^2$ and then ${\bar\lz}^{-1}=4/\pi^2\approx 0.405$.
 Thus, our estimates read as follows.
 $$\dz=0.25<\dz_1'={\dz_1^*}'=0.375< {\bar\lz}^{-1}\approx 0.405<\dz_1=\dz_1^*\approx 0.427<4\dz=1.$$

For this example, the results that $\dz_1=\dz_1^*$ and $\dz_1'={\dz_1^*}'$ are quite natural by symmetry. The not so obvious fact is $\dz_1\ge \dz_1^*$ in the most cases and the inequality can happen, as shown by numerical computations.

Actually, Proposition \ref{t11} is deduced from the next result which is an improvement of the main variational formula \rf{cmwf97}{Theorem 1}.

\thm\lb{t17} {\cms We have $\lz_1\ge {4}{\bar\lz}/ {D^2}$ and two variational formulas:
$$\aligned
{\bar\lz}&=
\sup_{f\in {\scr F}}\,\inf_{r\in (0, 1)}
\frac{f(r)}{\int_0^r C(s)^{-1}\d s \int_s^1 C(u)f(u)\d u}\\
&=\sup_{f\in {\scr F}}\,\inf_{r\in (0, 1)}
\frac{f(r)}{\int_r^1 C(s)\d s \int_0^s C(u)^{-1}f(u)\d u},
\endaligned$$
where  ${\scr F}=\big\{f\in {\scr C}[0, 1]: f|_{(0, 1)}>0\big\}$.
}
\dethm

Recall that
the estimates given in (\ref{08})---(\ref{10}) are all dimension-free. The next result
is an improvement of (\ref{09}) which depends on dimensions.

\prp\lb{t12}{\cms For $\lz_1$, we have lower bound (\ref{09}) replacing $K$ by $K M_{\az}$, where
 $$M_{\az} =\frac{\pi^2}{4} \int_0^1 (1-y)\cos\frac{\pi y}{2}\, \text{\rm sech}^2 (\az y) \d y$$
 regarding $\az$ as the constant $\az_0$ if $K>0$ and $|\az|\in (|\az_0|, \pi/2]$, where $|\az_0|$
(depending on $d$) is the first positive root of
\be\bigg(\frac{\pi}{2\sqrt{d-1}\,|\az|}+ \frac{\sqrt{d-1}\,|\az|}{2\pi}\bigg)\cos|\az|=1.\lb{15}\de}\deprp

In Proposition \ref{t12}, the improvement of (\ref{09}) is due to the fact
that $M_{\az}<1$ if $K<0$ and $M_{\az}>1$ if $K>0$. The proposition reduces to (\ref{09}) once $d\to\infty$
since then $\az\to 0$, and furthermore, $M_0=1$. Here is the picture of $M_{\az}$ (Figure 6).
\medskip

\begin{center}\includegraphics{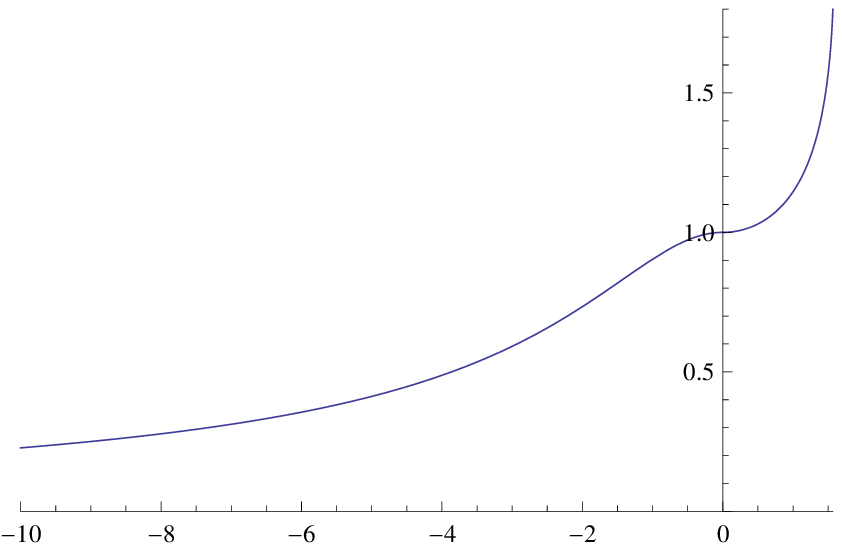}\newline
{{\bf Figure 6}\qd\rm The curve of $M_{\az}$ with $\az=\sqrt{-\text{sgn}(x)}\,|x|$, $x\in (-10, \pi/2)$. Obviously, $K<0$ iff $x<0$.}\end{center}
Note that the region on which Proposition \ref{t12} being available is smaller than that of (\ref{09}).
In view of (\ref{10}), this is reasonable since a larger lower bound can be held only
in a smaller region.

As mentioned in the earlier publication (cf. \rf{cmf05a}{Chapter 3} for instance)
or in Theorem \ref{t17}, our study on $\lz_1$ consists of two steps. The first one is reducing the higher dimensions to dimension one as shown by the first assertion of Proposition \ref{t11} or Theorem \ref{t17}. The second step is estimating ${\bar\lz}$. This was started by \ct{cmwf97}, continued by several papers mentioned before (\ref{08}),
and is also the main aim of the present paper to justify the power of our one-dimensional results \ct{cmf00}, \ct{cmf01} and \ct{cmf10}.
The lower bound $\dz_1^{-1}\vee {\dz_1^*}^{-1}$ of $\bar\lz$ provided by Proposition \ref{t11} is universal in the sense that the upper and lower bounds of $\bar\lz$ are the same up to a factor 4. Actually, all of a large number of examples we have ever computed,
as well as Figures 7--9 below, show that the ratio $\dz_1/\dz_1'$
(and $\dz_1^*/{\dz_1^*}'$) is no more than 2. It is somehow unexpected that the lower bound
$4 D^{-2}{\dz_1^*}^{-1}$ of $\lz_1$ is better than the others
except in two cases (cf. Figures 7--9 below). In the case that $K=0$, we have shown after Proposition \ref{t11}
by an example that $4 D^{-2}{\dz_1^*}^{-1}$ is not sharp. When $\az$ closes to $\pi/2$,
we are near the unit sphere and so $4 D^{-2}{\dz_1^*}^{-1}$ can not be better than (\ref{03}).
Thus, we may regard $4 D^{-2}{\dz_1^*}^{-1}$ (be careful to distinguish $\dz_1^*$ and $\dz_1$) as our general common lower bound, and regard (\ref{03}), (\ref{10}) and
Proposition \ref{t12} as an addition. We can now summarize the main result of the paper as follows.

\thm\lb{t13}{\cms In general, we have the following lower estimate:
$$\lz_1\!\ge\! \frac{4}{D^2}\bigg\{\frac{1}{\dz_1^*}\!\bigvee\! \sup_{s\in (0, 1)}\!s\bigg[{(1-s)\pi^2}\!- (d-1)\az^2 M_{\az}\bigg]\!\bigvee\! \bigg[ \mathbbold{1}_{\{K>0\}}\frac{d\, |\az|^2}{1-\cos^d |\az|}\bigg]\bigg\},$$
where $x\vee y=\max\{x, y\}$, $\dz_1^*$ and $M_{\az}$ are given in Propositions \ref{t11} and \ref{t12}, respectively, with a restriction on $|\az|$
for the middle term in the case of $K>0$: regarding
$\az$ as the constant $\az_0$ if $|\az|\in (|\az_0|, \pi/2]$, where $|\az_0|$ is the first positive root of (\ref{15}).
Besides, we also have the dimension-free lower bound (\ref{10}).
 The middle estimate is better than (\ref{10}) if
 $2\le d \le 7$, and conversely if $d\ge 10$.
}
\dethm

For the convenience of computation, according to (\ref{09}), we express the middle term in Theorem \ref{t13} as follows.
$$\aligned
&\sup_{s\in (0, 1)} s\big[{(1-s)\pi^2}\!- (d-1)\az^2 M_{\az}\big]\\
&\qd={
\begin{cases}
\big({\pi}/{2}-(d-1)\az^2 M_{\az}/(2\pi)\big)^2,\qd &-\pi^2\le -(d-1)\az^2 M_{\az}\le \pi^2\\
K M_{\az}, & -(d-1)\az^2 M_{\az}\in \big(\pi^2, (d-1) \pi^2/4\big]\\
0, &-(d-1)\az^2 M_{\az}<-\pi^2.
\end{cases}}
\endaligned$$
However, when $K>0$, $|\az|$ is essentially restricted to the subinterval $ (0, |\az_0|)$,
where $|\az_0|$ is the smallest positive root of equation (\ref{15}).

Fig. 7 illustrates the meaning of Theorem \ref{t13} ignoring the common factor $4D^{-2}$. Here, we consider only $K\ge 0$ and $d=2$.
Clearly, $\bar\lz$ is located between
the two dotted curves and the ratio of the upper and lower bounds \big(${\dz_1^*}^{\prime\,-1}$ and ${\dz_1^*}^{-1}$\big) of $\bar\lz$ is obviously less than 2. Note that here we use $*$ twice.
The curve, say Curve 1, corresponding to the last term in Theorem \ref{t13} is sharp at $\pi/2$ but is at the lowest level at origin; and the partially dashed curve, say Curve 2, corresponding to the middle term in Theorem \ref{t13} is sharp at $0$.
Note that Curve 2 is located above Curve 1 and is in particular higher than Curve 1 when $|\az|$ is close to $\pi/2$, which is impossible since Curve 1
is sharp at $\pi/2$ as we have just mentioned. Hence, a restriction on $|\az|$ for Curve 2 is necessary and what we adopted is $|\az|\le 0.97$ (that is ignoring the dashed part of Curve 2). Here for Curve 2, we omit the constant line starting from the endpoint of the non-dashed part of the curve to $\pi/2$ (cf. Fig. 9 below). It is interesting that Curves 1 and 2 together control the most part of the interval $[0, \pi/2]$. Usually, we have $\dz_1\ge \dz_1^*$. The exceptional is that here $\dz_1/ \dz_1^*\in (0.99993,\, 1)$ when $|\az|< 0.87$. Since this part is covered by Curve 2 and so $\dz_1$ is ignored in Theorem \ref{t13}.
\begin{center}
\includegraphics[width=10.0cm,height=6cm]{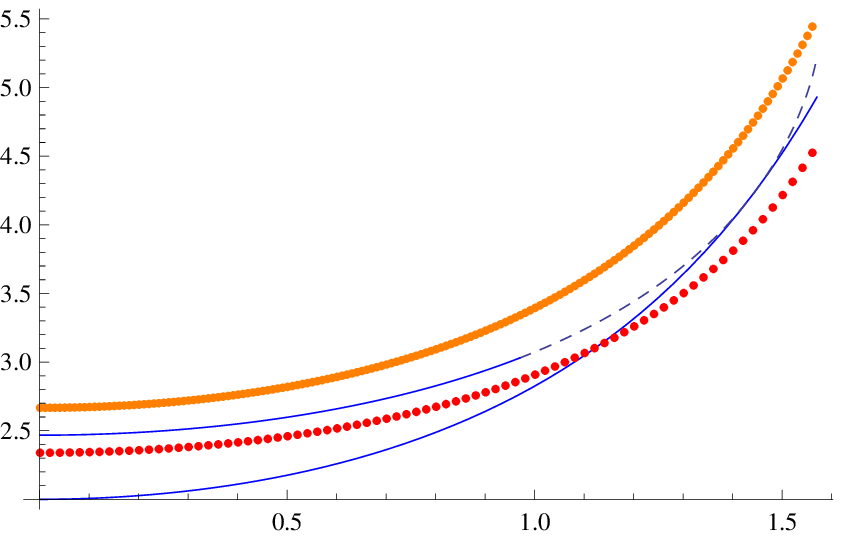}\newline
{{\bf Figure 7}\qd\rm The estimates of $\bar\lz$ in the case of $K\ge 0$ and $d=2$, $|\az|\le \pi/2$.}\end{center}

Fig. 8 represents the case of $K\le 0$ and $d=2$.
Again, the upper and lower bounds \big(${\dz_1'}^{-1}$ (but not ${\dz_1^*}^{\prime \, -1}$)
and ${\dz_1^*}^{-1}$\big) of $\bar\lz$ are given by the two of top dotted curves. The solid curve is determined by the middle term of Theorem \ref{t13}. The figure shows that the bounds ${\dz_1'}^{-1}$ and ${\dz_1^*}^{-1}$ are rather good,
even coincide each other for larger $\az$.
In a small interval around $0$, they are less sharper than the solid curve.
The dashed curve corresponds to (\ref{06}) and dotted part of the triangle corresponds to (\ref{07}). They are not comparable, and are less powerful than at least one of the others and so are disappeared in Theorem \ref{t13}. For $d>2$, the picture is similar but each curve decays fast (cf. Fig. 9).

Theorem \ref{t13} is stated unified in $K\in {\mathbb R}$. Fig. 9 is the result for both negative and positive $K$, where $\az=\sqrt{-\text{sgn}(x)}\,|x|$ with $x$ varies from $-2.5$ to $\pi/2$ and $d=5$. Note that the axes in these figures have different scales. The meaning of each curve should be clear. The top dotted curve is
${\dz_1'}^{-1}$ for $x<0$ and is ${\dz_1^*}^{\prime\,-1}$ for $x>0$. The part of the curves near $\pi/2$
is ignored, otherwise the left part of the curves would be very mixed. However, the shape
of the missed part is very much imaginable, up to 14.75 high.

\begin{center}
\includegraphics[width=10.0cm,height=6cm]{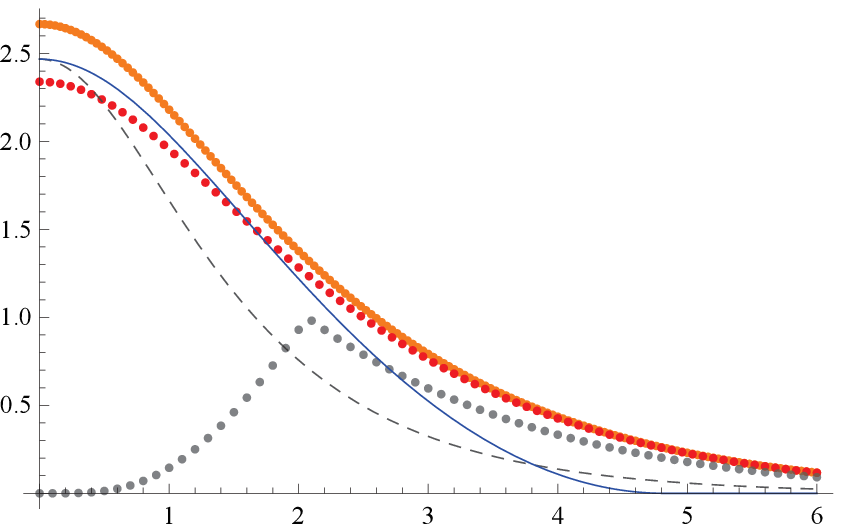}\newline
{{\bf Figure 8}\qd\rm The estimates of $\bar\lz$ in the case of $K\le 0$ and $d=2$, $0\le\az\le 6$.}\end{center}


\begin{center}
\includegraphics[width=10.0cm,height=8cm]{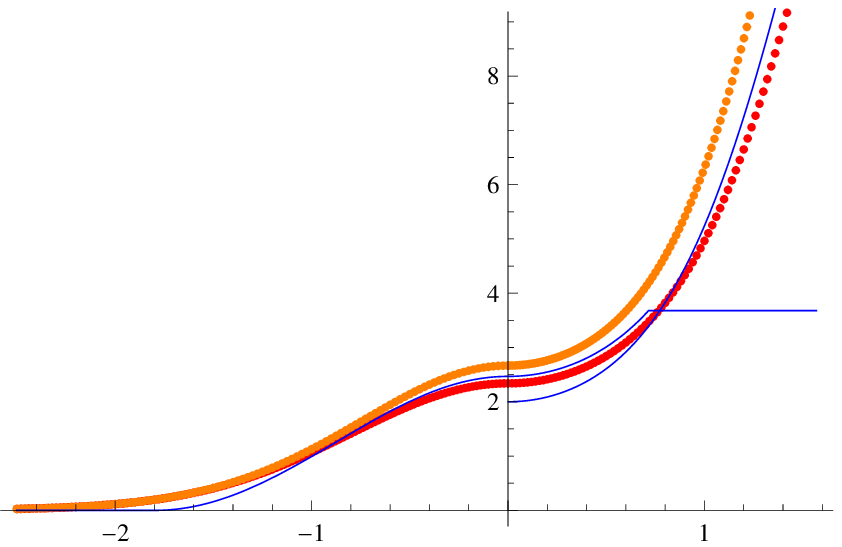}\newline
{{\bf Figure 9}\qd\rm The estimates of $\bar\lz$ when $d=5$ with
$\az=\sqrt{-\text{sgn}(x)}\,|x|$, $x\in (-2.5, \pi/2)$.}\end{center}

\rmk\lb{t14}{\rm Figures $7$--$9$ show that the ratio $\dz_1^*\big/{\dz_1^*}'\in [1, 4]$ (as well as ${{\dz_1^*}'}^{-1} - {\dz_1^*}^{-1}$) is controlled by its value at the endpoint $|\az|=\pi/2$. The ratio increases quickly from $1.2$ to $1.27$ and then slowly to $1.334$ ($<2$) when $d$ varies from $2$ to $5$ and then to $63$.}\dermk

\rmk\lb{t15}{\rm Consider the convex mean
$\ez_z= \gz_z {\dz_1^*}^{\prime\,-1} + \big(1-\gz_z\big){\dz_1^*}^{-1}$ with
$$\gz_0=\frac{5^{2/3}-5 \cdot 16^{-1}\pi^2}{5^{2/3}-10\cdot 3^{-1}}\approx 0.39\qqd\text{or}
\qqd \gz_{\pi/2}=\frac{4^{-1} d \pi^2-{\dz_1^*}^{-1}}{{\dz_1^*}^{\prime\,-1}-{\dz_1^*}^{-1}}\bigg|_{|\az|=\pi/2}.$$
The latter one depends on $d$ but the former one does not. Here $\ez_0$ is designed to be sharp $(=\pi^2/4)$ at $\az=0$ and so is $\ez_{\pi/2}(=d\, \pi^2/4)$ at $|\az|=\pi/2$ ($K> 0$).
In particular, $\gz_{\pi/2}\approx 0.35$ if $d=63.$ Numerical computations
(illustrating figures are given in the author's homepage) exhibit the following unexpected nice conclusion:
$${\bar\lz}-0.056\le \ez_{\pi/2}\le {\bar\lz}\le \ez_0\le {\bar\lz}+1.85,
\qqd 2\le d\le 63,\; \forall \az.$$
Thus, one may regard $\ez_0$ and $\ez_{\pi/2}$ (for each fixed $d$) as upper and lower bounds of ${\bar\lz}$, respectively, but they are almost the same since ${\bar\lz}\approx 155$ when $d=63$. This illustrates the power of Proposition $\ref{t11}$, and is independent of $(1)$--$(10)$.}\dermk

\section{Proofs}

\nnd{\bf Proof of Theorem \ref{t17}}. (a) Consider first the case that $\bar\lz$
is defined by its first equality given in the theorem. The first assertion is a comparison of $\lz_1$ with the principal eigenvalue $\lz_0$ of the operator
$${\overline L}=\frac{\d^2}{\d r^2} + (d-1)\az
\tanh(\az r)\frac{\d}{\d r}$$
on $(0, 1)$ with Dirichlet boundary at $0$ and Neumann boundary at $1$.
This was done in \ct{cmwf93},
as explicitly pointed out by \rf{cmwf97}{Remark d) after Theorem 1.1}.
Actually, the result came out by using the coupling method to a computation of some
distance which is certainly valued in the half-line. This reduces the higher dimensions to dimension one. Then it was proved in \rf{cmwf97}{Theorem 1.1}
that $\lz_0\ge \bar\lz$. The equality here holds because of \rf{cmf01}{Theorem 1.1}
noting that one allows the test function $f$ to be positive at origin not necessarily zero (cf. \rf{cmf01}{(1.3)}).

(b) Next, let $\lz_0^*$ be the principal eigenvalue of the dual operator
$${\overline L}^*=\frac{\d^2}{\d r^2} - (d-1)\az
\tanh(\az r)\frac{\d}{\d r}$$
with Neumann boundary at $0$ and Dirichlet boundary at $1$.
Here and in what follows, the notation ``$*$'' is used for dual quantity.
Then, we have not only $\lz_0=\lz_0^*$ but also the second equality for
$\bar\lz$ given in the theorem, as an analog of \rf{cmf10}{\S 5 and Theorem 2.4\,(3)}.
\deprf

\nnd{\bf Proof of Proposition \ref{t11}}. The first assertion comes from Theorem \ref{t17}. The result is
very helpful since the parameter $D$ is separated out, and then the three parameters $d$, $D$, and $K$ are now reduced to two:
$d$ and $\az$.

The ``$\dz$ part'' of the assertion was presented in
\rf{cmf05a}{Corollary 1.4}, comes originally from \rf{cmf00}{Theorem 2.2} with a change of the variable: $r\to r/D$ reducing the interval $(0, D)$ to $(0, 1)$. Actually, the original result is more general, including a vector field.
We mention that Proposition \ref{t11} is meaningful in such a general situation.

In view of the author's knowledge, the ``$\dz_1$ and $\dz_1'$ parts'' have not yet published in the geometric context. However, it is indeed
the first step of a general approximating procedure for the first eigenvalue, given by \rf{cmf01}{Theorem 2.2}. Here we state the procedure in the present context. Define $f_1=\sqrt{\fz}$,
$$f_{n+1}(r)=\int_0^r C(s)^{-1}\d s \int_s^1 C(u)f_n(u)\d u,\qqd n\ge 1$$
and $\dz_n=\sup_{r\in (0, 1)} f_{n+1}(r)/f_n(r)$. Then
$${\bar\lz}\ge \ldots\ge \dz_n^{-1}\ge \dz_{n-1}^{-1}\ge \ldots \ge \dz_1^{-1}\ge (4\dz)^{-1}.$$
Next, fix $r\in (0, 1)$ and define $f^{(r)}_1=\fz(\cdot\wedge r)$,
$$f^{(r)}_{n+1}=\int_0^{\bullet \wedge r} C(s)^{-1}\d s \int_s^1 C(u)f^{(r)}_n(u)\d u,\qqd n\ge 1$$
and $\dz_n'=\sup_{r\in (0, 1)}\inf_{s\in (0, 1)} f^{(r)}_{n+1}(s)\big/f^{(r)}_n(s)$. Then
$${\bar\lz}\le\ldots\le {\dz_n'}^{\!-1}\le {\dz_{n-1}'}^{\!\!\!\!\!\!\!\!-1}\le \ldots \le {\dz_1'}^{\!-1}\le \dz^{-1}.$$
Besides, we also have ${\bar\lz}\le {\bar\dz_n}^{-1}$ for all $n$, where
$${\bar\dz_n}= \sup_{r\in (0, 1)}
\int_0^1 f_n^{(r)}(s)^2 C(s)\d s\bigg/ \int_0^1 {\big(f_n^{(r)}(s)\big)'}^2 C(s)\d s,
\qqd n\ge 1.$$
Moreover, $\bar\dz_1=\dz_1'$. All of these approximating procedure comes from some
variational formulas \rf{cmf01}{Theorem 2.1}. In particular, the sequence $\{\dz_n\}_{n\ge 1}$ comes from the first variational
formula given in Theorem \ref{t17}.

The ``$\dz_1^*$ and ${\dz_1^*}'$ parts'' are parallel to
the ``$\dz_1$ and $\dz_1'$ parts'' based on the dual formula
just mentioned, as an analog of \rf{cmf10}{Theorem 3.2 and Corollary 3.3}.
The details are omitted here. We mention that the duality is studied more carefully in the author's forthcoming paper entitled ``Basic estimates of stability rate for one-dimensional diffusions''.

Clearly, the estimates given in Proposition \ref{t11} can be still improved by using the above approximating procedure.
\deprf

To prove Proposition \ref{t12}, we need some preparations. The following nice result it due to \rf{syzh}{proof of Lemma 2.2}, it is the key leading to (\ref{09}). As usual, we use ${\scr C}^{m}$ to denote the set
of functions having continuous $m$-th derivative.

\lmm\lb{t21} {\cms
Let $\lz$ and $f$ satisfy
$$f'' + F f'=-\lz f \qd \text{\cms on}\qd [0, \ell],\qqd \ell<\infty$$
with boundary conditions $f(0)=0$ and $f'(\ell)=0$, where
$F\in {\scr C}[0, \ell] \cap {\scr C}^1(0, \ell)$ having $F(0)=0$.
Then for each $s\in (0, 1)$, the function $g:=(f')^{2^{-1}(1-s)^{-1}}$
satisfy
$$4 s(1-s) \int_0^{\ell} {g'}^2 \d r = \int_0^{\ell} (\lz + s F') g^2 \d r.$$
Moreover, $g$ satisfies the mixed boundary condition: $g'(0)=0$ and $g(\ell)=0$.}
\delmm

The proof of Lemma \ref{t21} goes as follows. Making derivatives of the original equation, we get
$$f'''+ F f''= -({\lz} +F') f'.$$
Regarding $f'$ as a new function $h$, one sees that this corresponds to a second order Schr\"odinger operator (since the appearance of $F'$).
This step is known as an application of the coupling technique (cf. \rf{cmf10}{(10.4), (10.7)} and references therein). The next step we have used is to adopt the dual operator
${\d^2}/{\d x^2} -F {\d}/{\d x}$
(cf. \rf{cmf10}{(10.6)}) which is isospectral to the Schr\"odinger one. Refer to \rf{cmf10}{proof of Theorem 10.2} for more details. However, here the idea is different. Set $t=(1-s)^{-1}$.  Multiply the last
equation by $(f')^{t-1}$ and then integrate each term in the last equation over $(0, \ell)$. Based on the fact that $F(0)=0$ and then $f''(0)=0$, some careful computations lead to the required conclusion.

Having Lemma \ref{t21} at hand, the assertion (\ref{09}) is immediate. Simply replace $\lz + s F'$ by the constant $\lz + s\, \sup_{r\in (0, \ell)} F'(r)$ (which then does not depend on the dimension $d$ explicitly) and use the known exact inequality (cf. Proof of Proposition 1.2 below)
$$\int_0^{\ell}  g^2 \d r\le \bigg(\frac{2\ell}{\pi}\bigg)^2 \int_0^{\ell} {g'}^2 \d r. $$
It is at this step, the original proof is improved. This leads to
a use of the next result.

\lmm\lb{t22}{\cms
Let $a>0$ on $(0, 1)$ and $a\in {\scr C}^2(0, 1)$. Next, let $\hat\lz$ be the principal eigenvalue of the operator
${\widehat L}=a(x)^{-1}{\d^2}\big/{\d x^2}$
on $(0, 1)$ with Neumann boundary at $0$ and Dirichlet boundary at $1$.
Then $\hat\lz$ obeys the following estimates
$$\inf_{x\in (0, 1)}h(x)^{-1}\le \hat\lz \le \sup_{x\in (0, 1)}h(x)^{-1},$$
where
$$\aligned
h(x)
=\bigg(\frac{2}{\pi}\bigg)^2\bigg\{a(x)  -\sec\frac{\pi x}{2}
&\bigg[(1-x) a'(0)+(1-x)\int_0^x  a''(y) \cos\frac{\pi y}{2}\d y\\
&\;+\int_x^1 \big[-2 a'(y)+ (1-y) a''(y)\big] \cos\frac{\pi y}{2}\d y\bigg] \bigg\}.
\endaligned$$
}\delmm

\prf We are now in the case which is the continuous analog of \rf{cmf10}{\S 2}.
Noticing that when $a(x)$ is a constant, the assertion of the lemma is exact. In this case, the eigenfunction is cosine and so our approximation should start at the function cosine, because we are looking for such estimates which are sharp when $a(x)$ is a constant. Define $f_1(x)= \cos \frac{\pi x}{2}$ and
$$\aligned
f_n(x)&= \int_x^1 \d y \int_0^y a(u) f_{n-1}(u)\d u, \qqd x\in (0, 1),\; n\ge 1.
\endaligned$$
Here, the proof of the lemma is mainly for simplifying $f_2/f_1$ and so is rather
elementary. Exchanging the integrals, we get
$$\aligned
f_2(x)&=\int_x^1 \d y \int_0^y a(u) \cos\frac{\pi u}{2}\d u\\
&=\int_0^1 (1-x\vee u)  a(u) \cos\frac{\pi u}{2}\d u\\
&=(1-x)\int_0^x  a(y) \cos\frac{\pi y}{2}\d y
+\int_x^1  (1-y) a(y) \cos\frac{\pi y}{2}\d y.
\endaligned$$
As an application of the integration by parts formu1a, we have
$$\aligned
f_2(x)&= \frac{2}{\pi}\bigg[(1-x)\int_0^x  a(y)\, \d \sin\frac{\pi y}{2}
+\int_x^1  (1-y) a(y)\, \d \sin\frac{\pi y}{2}\bigg]\\
&= -\frac{2}{\pi}\bigg[
(1-x)\!\int_0^x\!\! a'(y) \sin\frac{\pi y}{2}\d y
 +\!\int_x^1\!\! \big[- a(y)+ (1-y) a'(y)\big] \sin\frac{\pi y}{2}\d y\bigg].
\endaligned$$
Using the integration by parts formu1a again, we get
$$\aligned
f_2(x)&=\bigg(\frac{2}{\pi}\bigg)^2\bigg\{a(x) \cos\frac{\pi x}{2}-
(1-x) a'(0) -(1-x)\int_0^x  a''(y) \cos\frac{\pi y}{2}\d y\\
&\qqd\qqd\qd -\int_x^1 \big[-2 a'(y)+ (1-y) a''(y)\big] \cos\frac{\pi y}{2}\d y \bigg\}.
\endaligned$$
Unlike the original one, no double integral appears in the last formula. Obviously,
${f_2(x)}/{f_1(x)}=h(x)$. The required assertion now follows by \rf{cmf01}{Theorem 1.1}.
\deprf

\nnd{\bf Proof of Proposition \ref{t12}}. Applying Lemma \ref{t21} to
${\overline L}$ and $\bar\lz$ defined in the proof of Theorem \ref{t17} with
$F (r)=(d-1)\,\az \;\text{tanh}(\az r),$
we obtain
\begin{align}
4 s(1-s) \int_0^{1} {g'}^2 \d r &= \int_0^{1}\big[\bar\lz + s (d-1) \az^2\text{sech}^2(\az r)\big] g^2 \d r\nnb\\
&={\bar\lz}\int_0^{1}  g^2 \d r+
s (d-1) \az^2 \int_0^{1}\text{sech}^2(\az r) g^2 \d r.\lb{11}
\end{align}

(a) Let $K\le 0$ and set $a (r)=\text{sech}^2(\az r)$. Then $a>0$ on $[0, 1]$ and
$$\aligned
a'(r)&=-2\,\az\, \text{sech}^2 (\az\, r)   \tanh (\az\, r),\\
a''(r)&=-2\,\az^2\, \text{sech}^4 (\az\, r)\,[2- \text{cosh} (2\az\, r)].
\endaligned$$
Applying Lemma \ref{t22} to this $a(r)$, we obtain, replacing $h$ by $h_{\az}$, that
$$\aligned
\frac{\pi^2}{4}h_{\az}(x)
=\text{sech}^2(\az x)&+2\az \sec \frac{\pi x}{2}
\bigg[\az (1-x)\int_0^x q_{\az}(y) \cos\frac{\pi y}{2}\d y\\
&\; + \int_x^1\!\big[-2 p_{\az}(y)+\az  (1-y) q_{\az}(y)\big] \cos\frac{\pi y}{2}\d y\bigg],
\endaligned$$
where
$$p_{\az}(y)\!=\text{sech}^2 (\az y)\, \text{tanh}(\az y)\qd\text{ and }\qd
q_{\az}(y)\!=\text{sech}^4 (\az y)\big[2-\text{cosh} (2\az y)\big].$$
Clearly, we have
\be \frac{\pi^2}{4}h_{\az}(0)=1+ 2\az \int_0^1 \big[-2 p_{\az}(y)+\az (1-y) q_{\az}(y) \big]\cos\frac{\pi y}{2}\d y.\lb{12} \de
Note that $\az^2>0$ iff $K<0$ in which case, $\sup_{x\in (0, 1)}h_{\az}(x)=h_{\az}(0)$
(see part (c) of the proof below).
Now, as an application of Lemma \ref{t22}, we have
$$\int_0^{1}\text{sech}^2(\az r) g^2 \d r
\le  h_{\az}(0) \int_0^{1} {g'}^2 \d r, \qqd K\le 0.$$
In particular, letting $\az\downarrow 0$, it follows that
$$\int_0^{1} g^2 \d r
\le  \frac{4}{\pi^2}\int_0^{1} {g'}^2 \d r.$$
By Lemma \ref{t22} again, we have for all $\az\ge 0$,
$$\az^2 \int_0^{1}\text{sech}^2(\az r) g^2 \d r
\le \az^2 h_{\az}(0) \int_0^{1} {g'}^2 \d r,$$
and then by (\ref{11}),
\be \pi^2 s(1-s)\le {\bar\lz}+s(d-1)\az^2\frac{\pi^2}{4} h_{\az}(0)
={\bar\lz}-s  \frac{K D^2}{4}\frac{\pi^2}{4} h_{\az}(0)\lb{13}.\de
Noticing that among
$q_{\az}$, $p_{\az}$, $\text{sech}^2(\az r)$, and
$\tanh (\az r)$, the former comes from the derivative of the latter,
starting from (\ref{12}), by using the integral by parts formula
three times, one finally sees that
$$\frac{\pi^2}{4} h_{\az}(0)=\frac{\pi^2}{4}
  \int_0^1 \bigg(\cos\frac{\pi y}{2} + \frac{\pi}{2}(1 - y)
  \sin \frac{\pi y}{2}\bigg) \frac{\tanh (\az y)}{\az} \d y.$$
The right-hand side gives us $M_{\az}$ in terms of the integration by parts
formula ($\sin=- \d \cos$). Applying this to (\ref{13}), Proposition \ref{t12} now follows from the first assertion of Proposition \ref{t11} in the case of $K\le 0$.

(b) Next, let $K>0$ and set
$$a(r)={\bar\lz}+s (d-1)\az^2\mbox{sech}^2(\az r)={\bar\lz}-s (d-1)|\az|^2\mbox{sec}^2(|\az| r).$$
In other words, we do not separate this $a$ into two parts as in (\ref{11}).
By (\ref{09}), in order that $a>0$ on $(0, 1)$, it suffices that
$$\bigg(\frac{\pi}{2\sqrt{d-1}\,|\az|}+ \frac{\sqrt{d-1}\,|\az|}{2\pi}\bigg)\cos |\az|>1.$$
Once the assertion is proved under this assumption, one may replace ``$>$'' here with ``$\ge$'' by
a limiting procedure. For the other part of $\az$, simply regard $\az$ as the constant
$\az_0$ (or $\az_0-\vz$ if necessary) since then $a(r)$ is upper bounded uniformly in $r$ by a positive constant on that subinterval of $\az$. We remark that even though the restriction here can be relaxed a little
we do not do so since the estimate is mainly essential for small $K$ (or for small $|\az|$). Applying Lemma \ref{t22} to this $a$, we obtain
$$\hat\lz^{-1}\le \sup_{x\in (0, 1)} h(x):\qd
h(x)=4{\bar\lz}/\pi^2 + s(d-1)\az^2 h_{\az}(x),$$
where $h_{\az}$ is the same as used in proof (a). Note that in the present case,
we have $\az^2<0$ and then
$$\sup_{x\in (0, 1)} h(x)=4{\bar\lz}/\pi^2+ s(d-1)\az^2 \inf_{x\in (0, 1)}h_{\az}(x).$$ This is the main different point to the previous case of $K\le 0$.
Luckily, we then have $\inf_{x\in (0, 1)}h_{\az}(x)=h_{\az}(0)$ (see part (c) of the proof below).
In this case, since $\az=i |\az|$, it may be more convenient to rewrite $h_{\az}$ as
$$\aligned
\frac{\pi^2}{4}h_{\az}(x)
=\text{sec}^2(|\az| x)&-2|\az| \sec \frac{\pi x}{2}
 \bigg[|\az| (1-x)\int_0^x q_{\az}^+(y) \cos\frac{\pi y}{2}\d y\\
&\; + \int_x^1\!\big[-2 p_{\az}^+(y)+|\az|  (1-y) q_{\az}^+(y)\big] \cos\frac{\pi y}{2}\d y\bigg],
\endaligned$$
where
$$p_{\az}^+(y)\!=\text{sec}^2 (|\az| y)\, \text{tan}(|\az| y)\qd\text{ and }\qd
q_{\az}^+(y)\!=\text{sec}^4 (|\az| y)\big[2-\text{cos} (2|\az| y)\big].$$
From this, one sees that $h_{\az}$ is always real for any $K\in {\mathbb R}$.
The remainder of the proof is similar to proof (a) above.

(c) To see that $\sup_{x\in (0, 1)}h_{\az}(x)=h_{\az}(0)$ when $K\le 0$ and
$\inf_{x\in (0, 1)}h_{\az}(x)=h_{\az}(0)$ when $K\ge 0$, it is helpful to
look at first three figures (Figures 10--12) for the latter case.
Fig. 10 shows that the surface is rather regular. But it may not be very clear that the minimum is attached at $x=0$ for each fixed $|\az|$, so two more figures (11 and 12) are included.
The pictures in case of $K\le 0$ are parallel. Based on the observation, it should not be hard to present an
analytic proof but we prefer to omit the details here for saving the space.
\deprf

Once again, as indicated in the proof of Lemma \ref{t22}, Proposition \ref{t12} uses only the first step
of our general approximating procedure for the first eigenvalue $\hat\lz$.
\medskip

\begin{center}\includegraphics{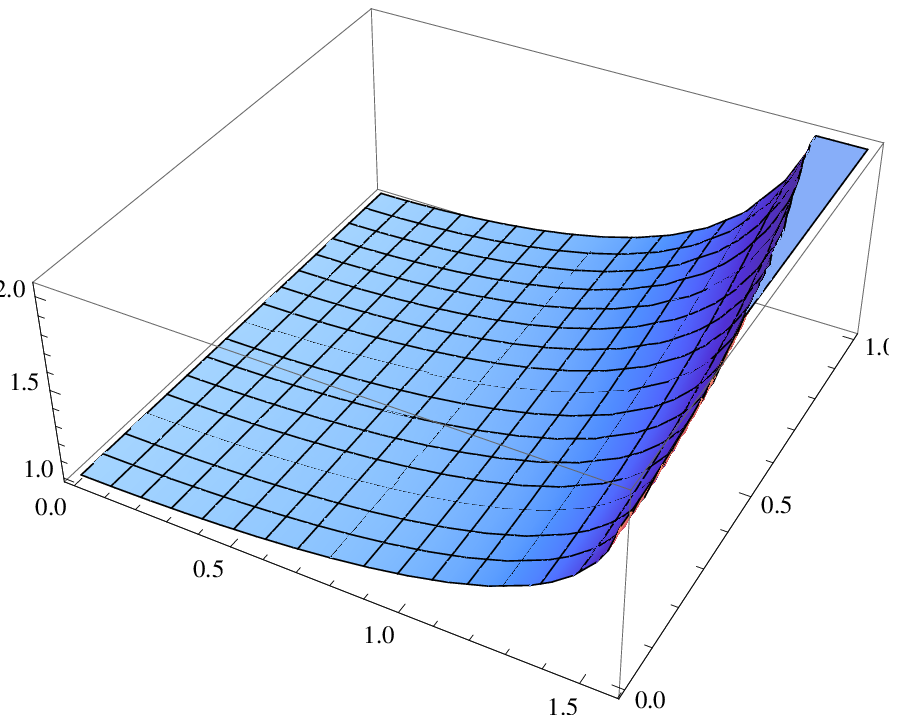}\newline
{{\bf Figure 10}\qd\rm When $K\ge 0$, the surface of $\pi^2 h_{\az}(x)/4$ in the interval $|\az|\in(0, 1.57075)$ and $x\in (0,1)$.}\end{center}

\begin{center}{{\includegraphics{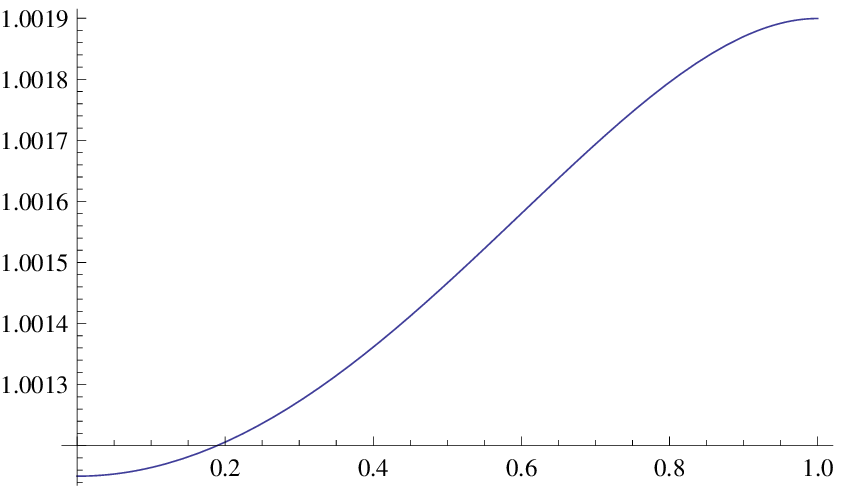}\hspace{-1.0truecm}}
\vskip-3.35truecm{
\hspace{-1.0truecm}{\includegraphics{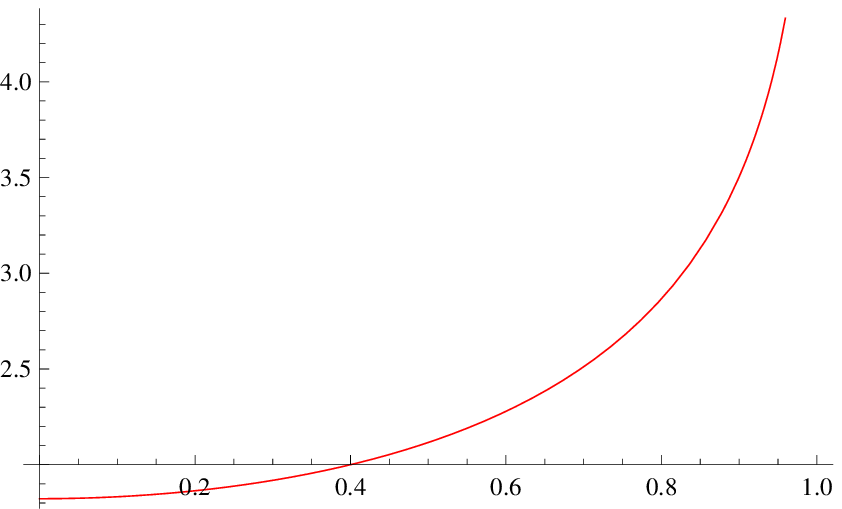}}
}\newline
{\bf Figure 11--12}\qd\rm The curve of $\pi^2 h_{\az}(x)/4$, $x \in (0, 1)$ when $K\ge 0$ for $|\az|=0.1$ or $|\az|=1.57$, respectively.}\end{center}

To conclude the paper, we make some remarks about the methods used in the paper.
Recall that all of the results (\ref{03}), (\ref{05})--(\ref{08}), and (\ref{10})
are an application of a coupling to a carefully designed (case by case) distance (equivalent to the Riemannian one). As mentioned in \rf{cmf94}{Theorem 6.2}, the method
works for more general ``cost'' functions, not necessarily a distance.
With coupling method in mind, the boundaries of the reduced process with operator ${\overline L}$ given in the proof of Theorem \ref{t17} is natural. Then the ``$\dz_1$''
(resp. ``$\dz_1'$'') part of Proposition \ref{t11} says that we do have a universal
distance which is equivalent to the Riemannian one and provides us a universal lower
bound $4 D^{-2} \dz_1^{-1}$. Thus, all of these results can be regarded as an application of the general variational formula given in \ct{cmwf97}. However, for part ``$\dz_1^*$'', as a dual of ${\overline L}$, it has a different probabilistic meaning.
The use of the dual technique is an essential new point of the present paper. Finally,
in Lemma \ref{t21} (or Proposition \ref{t12}), the original boundary conditions are also dualled, at the same time, the operator is changed. It is mainly a specific comparison result, one can not say that the two operators used in Lemma \ref{t21} have the same principal eigenvalue.

\medskip{\small

\nnd{\bf Acknowledgments}.
The main results of the paper have been presented at
``International Conference on Stochastic Partial Differential Equations and Related Topics (April 25 - 30, 2011; Chern Institute of Math.)'' and at
``Stochastic Analysis and Application to Financial Mathematics, in honor of Professor Jia-An Yan's 70th Birthday (July 4--6, 2011; AMSS, CAS)''. The author acknowledges the organizing groups for their kind invitation and financial support.}

\bigskip

\nnd {\small School of Mathematical Sciences, Beijing Normal University,\newline
Laboratory of Mathematics and Complex Systems (Beijing Normal University),\newline
\text{\qd} Ministry of Education, \newline
Beijing 100875, The People's Republic of China.\newline
E-mail: mfchen@bnu.edu.cn\newline
Home page:
    http://math.bnu.edu.cn/\~{}chenmf/main$\_$eng.htm
}

\end{document}